\documentclass[11pt]{article}
\usepackage{url}

\usepackage{a4,amsmath,amssymb}

\begin{document}

\title{Sums of generalized harmonic 
series for kids from five to fifteen}

\author{Zurab Silagadze\thanks{email: silagadze@inp.nsk.su} \\
{\small \em Budker Institute of Nuclear Physics SB RAS and} \\
{\small \em Novosibirsk State University, 630 090, Novosibirsk, Russia}}

\maketitle

\begin{abstract}
We examine the remarkable connection, first discovered by  Beukers, Kolk and
Calabi, between $\zeta(2n)$, the value of the Riemann zeta-function at an even
positive integer, and the volume of some $2n$-dimensional polytope. It can
be shown that this volume is equal to the trace of a compact self-adjoint
operator. We provide an explicit expression for the kernel of this operator
in terms of Euler polynomials. This explicit expression makes it easy
to calculate the volume of the polytope and hence $\zeta(2n)$. In the case of
odd positive integers, the expression for the  kernel enables us to rediscover
an integral representation for $\zeta(2n+1)$, obtained by a different method 
by Cvijovi\'{c} and Klinowski. Finally, we indicate that the origin of the 
miraculous Beukers-Kolk-Calabi change of variables in the multidimensional 
integral, which is at the heart of this circle of ideas, can be traced to 
the amoeba associated with the certain Laurent polynomial. The paper is 
dedicated to the memory of Vladimir Arnold (1937-2010). 
\end{abstract}

\section{Introduction}
In a nice little book \cite{1} Vladimir Arnold has collected 77 mathematical 
problems for kids from 5 to 15 to stimulate the development of a culture of
critical thinking in pupils. Problem 51 in this book asks the reader to 
calculate the sum of the inverse squares and prove Euler's celebrated formula
\begin{equation}
\sum\limits_{n=1}^\infty \frac{1}{n^2}=\frac{\pi^2}{6}.
\label{eq1}
\end{equation}
Well, there are many ways to do this (see, for example, \cite{2,3,4,5,5P,6,6A,
6B,7} and references therein), some maybe even accessible for kids under 
fifteen. However, in this note we concentrate on the approach of Beukers, 
Kolk and Calabi \cite{8}, further elaborated by Elkies in \cite{9}. This 
approach incorporates pleasant features which all the kids (and even some 
adults) adore: simplicity, magic and the depth that allows one to go beyond 
the particular case (\ref{eq1}). The simplicity, however, is not everywhere 
explicit in \cite{8} and \cite{9}, while the magic longs for explanation 
after the first admiration fades away. Below we will try to enhance the 
simplicity of the approach and somewhat uncover the secret of the magic.

The paper is organized as follows. In the first two sections we reconsider the
evaluation of $\zeta(2)$ and $\zeta(3)$ so that technical details of the 
general case do not obscure the simple underlying ideas. Then we elaborate 
the general case and give the main result of this work, the formula for 
the  kernel which allows us to simplify considerably the evaluation of 
$\zeta(2n)$ from \cite{8,9} and re-derive Cvijovi\'{c} and Klinowski's 
integral representation \cite{10} for $\zeta(2n+1)$. Finally, we ponder over 
the mysterious relations between the sums of generalized harmonic series and 
amoebas, first indicated by Passare in \cite{7}. This relation enables us to 
uncover somewhat the origin of the Beukers-Kolk-Calabi's highly non-trivial 
change of variables.

\section{Evaluation of $\zeta(2)$}
Recall the definition of the Riemann zeta function
\begin{equation}
\zeta(s)=\sum\limits_{n=1}^\infty \frac{1}{n^s}.
\label{eq2}
\end{equation}
The sum (\ref{eq1}) is just $\zeta(2)$ which we will now evaluate
following the method of Beukers, Kolk and Calabi \cite{8}. Our starting point
will be the dilogarithm function
\begin{equation}
Li_2(x)=\sum\limits_{n=1}^\infty \frac{x^n}{n^2}.
\label{eq3}
\end{equation}
Clearly, $Li_2(0)=0$ and $Li_2(1)=\zeta(2)$. Differentiating (\ref{eq3}), we
get
$$x\,\frac{d}{dx}Li_2(x)=\sum\limits_{n=1}^\infty \frac{x^n}{n}=
-\ln{(1-x)},$$
and, therefore,
\begin{equation}
\zeta(2)=Li_2(1)=-\int\limits_0^1\frac{\ln{(1-x)}}{x}dx=
\iint\limits_\Box\frac{dx\,dy}{1-xy},
\label{eq4}
\end{equation}
where $\Box=\{(x,y): 0\le x\le 1, 0\le y \le 1\}$ is the unit square.
Let us note 
\begin{equation}
\iint\limits_\Box\frac{dx\,dy}{1-xy}+
\iint\limits_\Box\frac{dx\,dy}{1+xy}=
2\iint\limits_\Box\frac{dx\,dy}{1-x^2y^2},
\label{eq5}
\end{equation}
and
\begin{equation}
\iint\limits_\Box\frac{dx\,dy}{1-xy}-
\iint\limits_\Box\frac{dx\,dy}{1+xy}=
\frac{1}{2}\iint\limits_\Box\frac{dx\,dy}{1-xy},
\label{eq6}
\end{equation}
where the last equation follows from
$$\iint\limits_\Box\frac{2xy}{1-x^2y^2}dx\,dy=
\frac{1}{2}\iint\limits_\Box\frac{d(x^2)\,d(y^2)}{1-x^2y^2}=
\frac{1}{2}\iint\limits_\Box\frac{dx\,dy}{1-xy}.$$
It follows from equations (\ref{eq5}) and (\ref{eq6}) that
\begin{equation}
\zeta(2)=\frac{4}{3}\iint\limits_\Box\frac{dx\,dy}{1-x^2y^2}.
\label{eq7}
\end{equation}
Now let us make the magic Beukers-Kolk-Calabi change of variables in this
two-dimensional integral \cite{8}
\begin{equation}
x=\frac{\sin{u}}{\cos{v}},\;\;\; y=\frac{\sin{v}}{\cos{u}},
\label{eq8}
\end{equation}
with Jacobian determinant
$$\frac{\partial(x,y)}{\partial(u,v)}=\left | \begin{array}{cc} 
\frac{\cos{u}}{\cos{v}} &
\frac{\sin{v}\,\sin{u}}{\cos^2{u}} \\ & \\
\frac{\sin{u}\,\sin{v}}{\cos^2{v}} & \frac{\cos{v}}{\cos{u}} 
\end{array} \right |=1-\frac{\sin^2{u}\sin^2{v}}{\cos^2{v}\cos^2{u}}=
1-x^2y^2.$$
Then miraculously
\begin{equation}
\zeta(2)=\frac{4}{3}\iint \limits_\Delta du\,dv=\frac{4}{3}\,Area(\Delta),
\label{eq9}
\end{equation}
where $\Delta$ is the image of the unit square $\Box$ under the transformation
$(x,y)\to (u,v)$. It is easy to show that $\Delta$ is the  isosceles right 
triangle $\Delta=\{(u,v): u\ge 0, v\ge 0, u+v\le\pi/2\}$ and, therefore,
\begin{equation}
\zeta(2)=\frac{4}{3}\,\frac{1}{2}\left (\frac{\pi}{2}\right )^2=
\frac{\pi^2}{6}.
\label{eq10}
\end{equation}
``Beautiful -- even more so, as the same method of proof extends to the
computation of $\zeta(2k)$ in terms of a $2k$-dimensional integral, 
for all $k\ge 1$'' \cite{11}. However, before considering the general case,
we check whether the trick works for $\zeta(3)$.

\section{Evaluation of $\zeta(3)$}  
In the case of $\zeta(3)$, we begin with trilogarithm
\begin{equation}
Li_3(x)=\sum\limits_{n=1}^\infty \frac{x^n}{n^3},
\label{eq11}
\end{equation}
and using
$$x\,\frac{d}{dx}Li_3(x)=Li_2(x)=-\int\limits_0^x\frac{\ln{(1-y)}}{y}dy,$$
we get
\begin{equation}
\zeta(3)=Li_3(1)=-\int\limits_0^1\frac{dx}{x}\int\limits_0^x
\frac{\ln{(1-y)}}{y}dy.
\label{eq12}
\end{equation}
But
$$-\frac{1}{x}\int\limits_0^x\frac{\ln{(1-y)}}{y}dy=-\int\limits_0^1
\frac{\ln{(1-xz)}}{xz}dz=\int\limits_0^1 dz \int\limits_0^1 \frac{dy}
{1-xyz},$$
and finally
\begin{equation}
\zeta(3)=Li_3(1)=\iiint\limits_{\Box_3}\frac{dx\,dy\,dz}{1-xyz},
\label{eq13}
\end{equation}
where $\Box_3=\{(x,y,z): 0\le x \le 1, 0\le y \le 1, 0\le z \le 1\}$ is the
unit cube. By a similar trick as before, we can transform (\ref{eq13})
into the integral
\begin{equation}
\zeta(3)=\frac{8}{7}\iiint\limits_{\Box_3}\frac{dx\,dy\,dz}{1-x^2y^2z^2},
\label{eq14}
\end{equation}
and here the analogy with the previous case ends, unfortunately, because the
generalization of the Beukers-Kolk-Calabi change of variables does not lead
in this case to a simple integral. However, it is interesting to note that
the hyperbolic version of this change of variables
\begin{equation}
x=\frac{\sinh{u}}{\cosh{v}},\;\;\;y=\frac{\sinh{v}}{\cosh{w}},\;\;\;
z=\frac{\sinh{w}}{\cosh{u}}
\label{eq15}
\end{equation}
{\it does} indeed produce an interesting result
\begin{equation}
\zeta(3)=\frac{8}{7}\iiint\limits_{U_3} du\,dv\,dw = \frac{8}{7}\, 
\mathrm{Vol}(U_3),
\label{eq16}
\end{equation}
where $U_3$ is a complicated 3-dimensional shape defined by the inequalities
$$ u\ge 0, \;\;v\ge 0,\;\; w\ge 0,\;\; \sinh{u}\le \cosh{v},\;\;
\sinh{v}\le \cosh{w}, \;\; \sinh{w}\le \cosh{u}.$$
Unfortunately, unlike the previous case, there is no obvious simple way  
to calculate the volume of $U_3$.

However, there is a second way to convert the integral (\ref{eq12}) for
$\zeta(3)$ in which the Beukers-Kolk-Calabi change of variables still plays
a helpful role. We begin with the identity
\begin{equation}
\zeta(3)=-\int\limits_0^1\frac{dx}{x}\int\limits_0^x
\frac{\ln{(1-y)}}{y}dy=-\iint\limits_D\frac{\ln{(1-y)}}{xy}\,dx\,dy,
\label{eq17}
\end{equation}
where the domain of the 2-dimensional integration is the triangle $D=\{(x,y): 
x\ge 0, y\ge 0, y\le x\}$. Interchanging the order of  
integration in (\ref{eq17}), we get
$$\zeta(3)=-\int\limits_0^1\frac{\ln{(1-y)}}{y}\,dy\int\limits_y^1\frac{dx}
{x},$$
which can be transformed further as follows
$$\zeta(3)=\int\limits_0^1\frac{\ln{(1-y)}\ln{y}}{y}\,dy=-\int\limits_0^1
\ln{y}\,dy\int\limits_0^1\frac{dx}{1-xy}=-\iint\limits_\Box \frac{\ln{y}}
{1-xy}\,dx\,dy,$$
or in a more symmetrical form
\begin{equation}
\zeta(3)=-\frac{1}{2}\iint\limits_\Box \frac{\ln{(xy)}}
{1-xy}\,dx\,dy.
\label{eq18}
\end{equation}
Note that
$$\iint\limits_\Box \frac{2xy\,\ln{(xy)}}{1-x^2y^2}\,dx\,dy=
\frac{1}{4}\iint\limits_\Box \frac{\ln{(x^2y^2)}}{1-x^2y^2}\,d(x^2)\,d(y^2)=
\frac{1}{4}\iint\limits_\Box \frac{\ln{(xy)}}{1-xy}\,dx\,dy.$$
Therefore, we can modify (\ref{eq5}) and (\ref{eq6}) accordingly and using
them transform (\ref{eq18}) into
\begin{equation}
\zeta(3)=-\frac{4}{7}\iint\limits_\Box \frac{\ln{(xy)}}{1-x^2y^2}\,dx\,dy.
\label{eq19}
\end{equation}
At this point we can use the Beukers-Kolk-Calabi change of variables 
(\ref{eq8}) in (\ref{eq19}) and as a result we get 
\begin{equation}
\zeta(3)=-\frac{4}{7}\iint\limits_\Delta\ln{(\tan{u}\,\tan{v})}\,du\,dv=
-\frac{8}{7}\iint\limits_\Delta\ln{(\tan{u})}\,du\,dv.
\label{eq20}
\end{equation}
But this equation indicates that
$$\zeta(3)=-\frac{8}{7}\int\limits_0^{\pi/2}du\,\ln{(\tan{u})}\int\limits_
0^{\pi/2-u} dv=-\frac{8}{7}\int\limits_0^{\pi/2}\left (\frac{\pi}{2}-u
\right )\ln{(\tan{u})}\,du,$$
which after the substitution $x=\frac{\pi}{2}-u$ becomes
\begin{equation}
\zeta(3)=-\frac{8}{7}\int\limits_0^{\pi/2}x\ln{(\cot{x})}\,dx=
\frac{8}{7}\int\limits_0^{\pi/2}x\ln{(\tan{x})}\,dx.
\label{eq21}
\end{equation}
But
$$\int\limits_0^{\pi/2}\ln{(\tan{x})}\,dx=-\int\limits_{\pi/2}^0
\ln{(\cot{u})}\,du=-\int\limits_0^{\pi/2}\ln{(\tan{u})}\,du=0,$$
which allows us to rewrite (\ref{eq21}) as follows
$$\zeta(3)=\frac{8}{7}\int\limits_0^{\pi/2}\left (x-\frac{\pi}{4}\right )
\ln{(\tan{x})}\,dx=\frac{8}{7}\int\limits_0^{\pi/2}\ln{(\tan{x})}\,
\frac{d}{dx}\left (\frac{x^2}{2}-\frac{\pi}{4}x\right )\,dx,$$
and after integration by parts and rescaling $x\to x/2$ we end with
\begin{equation}
\zeta(3)=\frac{1}{7}\int\limits_0^\pi \frac{x(\pi-x)}{\sin{x}}\,dx.
\label{eq22}
\end{equation}
This is certainly an interesting result. Note that  until quite recently very 
few definite integrals of this kind, involving cosecant or secant functions, 
were known and present in standard tables of integrals \cite{12,13,14}. In 
fact (\ref{eq22}) is a special case of the more general result \cite{10} 
which we are going now to establish.

\section{The general case of $\zeta(2n)$}
The evaluation of $\zeta(2)$ can be straightforwardly generalized. The 
polylogarithm function
\begin{equation}
Li_s(x)=\sum\limits_{n=1}^\infty \frac{x^n}{n^s}
\label{eq23}
\end{equation}
obeys
$$x\,\frac{d}{dx}Li_s(x)=Li_{s-1}(x),$$
and hence
\begin{equation}
Li_s(x)=\int\limits_{0}^x \frac{Li_{s-1}(y)}{y}\,dy.
\label{eq24}
\end{equation}
Repeated application of this identity allows to write
\begin{equation}
\zeta(n)=Li_n(1)=\int\limits_0^1\frac{dx_1}{x_1}\int\limits_0^{x_1}
\frac{dx_2}{x_2}\ldots\int\limits_0^{x_{n-2}}\frac{dx_{n-1}}{x_{n-1}}
\left [-\ln{(1-x_{n-1})}\right ].
\label{eq25}
\end{equation}
After rescaling
$$x_1=y_1,\;\;x_2=x_1y_2,\;\;x_3=x_2y_3,\ldots, x_{n-1}=x_{n-2}y_{n-1}=
y_1y_2\cdots y_{n-1},$$
and using
$$\int\limits_0^1\frac{dy_n}{1-y_1y_2\cdots y_n}=-\frac{1}{y_1y_2\cdots 
y_{n-1}}\,\ln{(1-y_1y_2\cdots y_{n-1})},$$
we get
\begin{equation}
\zeta(n)=\idotsint\limits_{\Box_n}\frac{dy_1\,dy_2\cdots dy_n}
{1-y_1y_2\cdots y_n},
\label{eq26}
\end{equation}
where $\Box_n$ is $n$-dimensional unit hypercube. The analogs of (\ref{eq5})
and (\ref{eq6}) are
$$\idotsint\limits_{\Box_n}\frac{dx_1\cdots dx_n}
{1-x_1\cdots x_n}+\idotsint\limits_{\Box_n}\frac{dx_1\cdots dx_n}
{1+x_1\cdots x_n}=2\idotsint\limits_{\Box_n}\frac{dx_1\cdots dx_n}
{1-x^2_1\cdots x^2_n}$$
and
$$\idotsint\limits_{\Box_n}\frac{dx_1\cdots dx_n}
{1-x_1\cdots x_n}-\idotsint\limits_{\Box_n}\frac{dx_1\cdots dx_n}
{1+x_1\cdots x_n}=\frac{1}{2^{n-1}}\idotsint\limits
_{\Box_n}\frac{dx_1\cdots dx_n}{1-x_1\cdots x_n},$$
from which it follows that (\ref{eq26}) is equivalent to
\begin{equation}
\zeta(n)=\frac{2^n-1}{2^n}\idotsint\limits_{\Box_n}\frac{dx_1\cdots dx_n}
{1-x^2_1\cdots x^2_n}.
\label{eq27}
\end{equation}
If we now make a change of variables that generalizes (\ref{eq8}), namely
\begin{equation}
x_1=\frac{\sin{u_1}}{\cos{u_2}},\;\;x_2=\frac{\sin{u_2}}{\cos{u_3}},\ldots,
\;x_{n-1}=\frac{\sin{u_{n-1}}}{\cos{u_n}},\;\;x_n=\frac{\sin{u_n}}{\cos{u_1}}.
\label{eq28}
\end{equation}
we, in general, encounter a problem because the Jacobian of (\ref{eq28}) is
\cite{8,9}
$$\frac{\partial(x_1,\ldots,x_n)}{\partial(u_1,\ldots,u_n)}=
1-(-1)^n\,x^2_1x^2_2\cdots x^2_n,$$
and, therefore, only for even $n$ we will get a ``simple'' integral. For the
hyperbolic version of (\ref{eq28}),
\begin{equation}
x_1=\frac{\sinh{v_1}}{\cosh{v_2}},\;\;x_2=\frac{\sinh{v_2}}{\cosh{v_3}},
\ldots, x_{n-1}=\frac{\sinh{v_{n-1}}}{\cosh{v_n}},\;\;x_n=\frac{\sinh{v_n}}
{\cosh{v_1}},
\label{eq29}
\end{equation}
the Jacobian has the ``right'' form
$$\frac{\partial(x_1,\ldots,x_n)}{\partial(v_1,\ldots,v_n)}=
1-x^2_1x^2_2\cdots x^2_n,$$
and we get
\begin{equation}
\zeta(n)=\frac{2^n}{2^n-1}\idotsint\limits_{U_n}dv_1\cdots dv_n=
\frac{2^n}{2^n-1}\;\mathrm{Vol}_{n}(U_n).
\label{eq30}
\end{equation}
However, the figure $U_n$ has a complicated shape and it is not altogether
clear how to calculate its $n$-dimensional volume $\mathrm{Vol}_{n}(U_n)$
(nevertheless, a hyperbolic version can lead to some new insights 
\cite{6A,6B}). Therefore, for a moment, we concentrate on the even values 
of $n$ for which (\ref{eq28}) works perfectly well and leads to \cite{8,9} 
\begin{equation}
\zeta(2n)=\frac{2^{2n}}{2^{2n}-1}\idotsint\limits_{\Delta_{2n}}du_1\cdots 
du_n=\frac{2^{2n}}{2^{2n}-1}\;\mathrm{Vol}_{2n}(\Delta_{2n}),
\label{eq31}
\end{equation}
where $\Delta_n$ is a $n$-dimensional polytope defined through the 
inequalities
\begin{equation}
\Delta_n=\left \{(u_1,\ldots,u_n): u_i\ge 0,\; u_i+u_{i+1}\le \frac{\pi}{2}
\right \}.
\label{eq32}
\end{equation}
It is assumed in (\ref{eq32}) that $u_i$ are indexed cyclically (mod $n$) and
therefore $u_{n+1}=u_1$.

There exists an elegant method due to Elkies \cite{9} for calculating the 
$n$-volume of $\Delta_n$ (earlier calculations of this type can be found in
\cite{15}). Obviously
\begin{equation}
\mathrm{Vol}_{n}(\Delta_{n})=\left (\frac{\pi}{2}\right )^n\,
\mathrm{Vol}_{n}(\delta_{n}),
\label{eq33}
\end{equation}
where $\mathrm{Vol}_{n}(\delta_{n})$ is the $n$-dimensional volume of the 
rescaled polytope
\begin{equation}
\delta_n=\left \{(u_1,\ldots,u_n): u_i\ge 0,\; u_i+u_{i+1}\le 1
\right \}.
\label{eq34}
\end{equation}
If we introduce the characteristic function $K_1(u,v)$ of the isosceles right 
triangle $\{(u,v): u,v\ge 0,\,u+v\le 1\}$ that is 1 inside  the triangle and 
0 outside of it, then \cite{9}
$$\mathrm{Vol}_{n}(\delta_{n})=\int\limits_0^1\ldots\int\limits_0^1
\prod\limits_{i=1}^nK_1(u_i,u_{i+1})\,du_1\ldots du_n=\int\limits_0^1du_1
\int\limits_0^1du_2 \,K_1(u_1,u_2)\ldots$$
\begin{equation}
\int\limits_0^1du_{n-1}\,K_1(u_{n-2},u_{n-1})\int\limits_0^1du_n\,
K_1(u_{n-1},u_n)\,K_1(u_n,u_1).
\label{eq35}
\end{equation}
Let us note that $K_1(u,v)$ can be interpreted \cite{9} as the kernel of the
linear operator $\hat T$ on the Hilbert space $L^2(0,1)$, defined as follows
\begin{equation}
(\hat T f)(u)=\int\limits_0^1 K_1(u,v)f(v)\,dv=\int\limits_0^{1-u} f(v)\,dv.
\label{eq36}
\end{equation}
Then (\ref{eq35}) shows that $\mathrm{Vol}_{n}(\delta_{n})$ equals just to 
the trace of the operator $\hat T^n$:
\begin{equation}
\mathrm{Vol}_{n}(\delta_{n})=\int\limits_0^1 K_n(u_1,u_1)\,du_1,
\label{eq37}
\end{equation}
whose kernel $K_n(u,v)$ obeys the recurrence relation
\begin{equation}
K_n(u,v)=\int\limits_0^1K_1(u,u_1)\,K_{n-1}(u_1,v)\,du_1.
\label{eq38}
\end{equation}
Surprisingly, we can find a simple enough solution of this recurrence 
relation \cite{15S}. Namely,
\begin{eqnarray} && \hspace{15mm}
K_{2n}(u,v)=(-1)^n\,\frac{2^{2n-2}}{(2n-1)!} \times  \nonumber \\ &&
\left\{ \left [  E_{2n-1}\left (
\frac{u+v}{2}\right )+E_{2n-1}\left (\frac{u-v}{2}\right )\right ]
\theta(u-v)+\right . \nonumber \\ && \left .  
\hspace*{2.5mm}\left [  E_{2n-1}\left (\frac{u+v}{2}\right )+
E_{2n-1}\left (\frac{v-u}{2}\right )\right ]\theta(v-u)\right\},
\label{eq39}
\end{eqnarray}
and
\begin{eqnarray} && \hspace{15mm}
K_{2n+1}(u,v)=(-1)^n\,\frac{2^{2n-1}}{(2n)!} \times  \nonumber \\ &&
\left\{ \left [  E_{2n}\left (
\frac{1-u+v}{2}\right )+E_{2n}\left (\frac{1-u-v}{2}\right )\right ]
\theta(1-u-v)+\right . \nonumber \\ && \left .  
\hspace*{2.5mm}\left [  E_{2n}\left (\frac{1-u+v}{2}\right )-
E_{2n}\left (\frac{u+v-1}{2}\right )\right ]\theta(u+v-1)\right\}.
\label{eq40}
\end{eqnarray}
In these formulas $E_n(x)$ are the Euler polynomials \cite{16} and $\theta(x)$ 
is the Heaviside step function
$$\theta(x)=\left\{\begin{array}{c} 1,\;\; {\mathrm if}\;\; x>0, \\ \\
\frac{1}{2},\;\; {\mathrm if}\;\; x=0, \\ \\ 0,\;\; {\mathrm if}\;\; x<0. 
\end{array} \right .$$
After they are guessed, it is quite straightforward to prove (\ref{eq39}) and
(\ref{eq40}) by induction using the recurrence relation (\ref{eq38}) and the 
following properties of the Euler polynomials:
\begin{equation}
\frac{d}{dx}E_n(x)=nE_{n-1}(x),\;\;\; E_n(1-x)=(-1)^nE_n(x).
\label{eq41}
\end{equation}
In particular, after rather lengthy but straightforward integration we get
$$\int\limits_0^{1-u}K_{2n+1}(u_1,v)du_1=K_{2n+2}(u,v)-X,$$
where
$$X=(-1)^{n+1}\,\frac{2^{2n}}{(2n+1)!}
\left[ E_{2n+1}\left (\frac{1+v}{2}\right )+E_{2n+1}\left (\frac{1-v}
{2}\right )\right ].$$
But 
$$\frac{1-v}{2}=1-\frac{1+v}{2}$$
and the second identity of (\ref{eq41}) then implies that $X=0$. 

Therefore the only relevant question is how (\ref{eq39}) and (\ref{eq40}) 
were guessed. Maybe the best way to explain the ``method'' used is to refer 
to problem 13 from the aforementioned book \cite{1}. To demonstrate the 
cardinal difference  between the ways problems are posed and solved by 
physicists and by mathematicians, Arnold provides the following problem for 
children:

``On a bookshelf there are two volumes of Pushkin's poetry. The thickness 
of the pages of each volume is 2 cm and that of each cover 2 mm. A worm bores 
through from the first page of the first volume to the last page of the 
second, along the normal direction to the pages. What distance did it cover?''

Usually kids have no problems to find the unexpected correct answer, 4 mm,
in contrast to adults. For example, the editors of the highly respectable 
physics journal initially corrected the text of the problem itself into: 
``from the last page of first volume to the first page of the second'' to 
``match'' the answer given by Arnold \cite{1,17}. The secret of kids lies in 
the experimental method used by them: they simple go to the shelf and see how 
the first page of the first volume and the last page of the second are 
situated with respect to each other.

The method that led to (\ref{eq39}) and (\ref{eq40}) was exactly of this 
kind: we simply calculated a number of explicit expressions for $K_n(u,v)$
using (\ref{eq38}) and tried to locate regularities in these expressions.

Having (\ref{eq39}) at our disposal, it is easy to calculate the integral in
(\ref{eq37}). Namely, because
\begin{equation}
K_{2n}(u,u)=(-1)^n\,\frac{2^{2n-2}}{(2n-1)!}\left [ E_{2n-1}(u)+ E_{2n-1}(0)
\right ],
\label{eq42}
\end{equation}
and
\begin{equation}
E_{2n-1}(u)=\frac{1}{2n}\,\frac{d}{du}\,E_{2n}(u),
\label{eq43}
\end{equation}
we get 
\begin{equation}
\mathrm{Vol}_{2n}(\delta_{2n})=\int\limits_0^1 K_{2n}(u,u)\,du=
(-1)^n\,\frac{2^{2n-2}}{(2n-1)!}E_{2n-1}(0),
\label{eq44}
\end{equation}
(note that $E_{2n}(0)=E_{2n}(1)=0$.)
But $E_{2n-1}(0)$ can be expressed in terms of the Bernoulli numbers 
\begin{equation}
E_{2n-1}(0)=-\frac{2}{2n}\,(2^{2n}-1)B_{2n},
\label{eq45}
\end{equation}
and combining (\ref{eq31}), (\ref{eq33}), (\ref{eq44}) and (\ref{eq45}), we 
finally reproduce the celebrated formula
\begin{equation}
\zeta(2n)=(-1)^{n+1}\,\frac{2^{2n-1}}{(2n)!}\,\pi^{2n}\,B_{2n}.
\label{eq46}
\end{equation}

\section{The general case of $\zeta(2n+1)$}
The evaluation of $\zeta(3)$ can be also generalized straightforwardly. We 
have
$$\zeta(n)=\int\limits_0^1\frac{Li_{n-1}(x_1)}{x_1}dx_1=
\int\limits_0^1\frac{dx_1}{x_1}\int\limits_0^{x_1}
\frac{Li_{n-2}(x_2)}{x_2}\,dx_2=\iint\limits_D\frac{Li_{n-2}(x_2)}
{x_1x_2}dx_1dx_2.$$
Interchanging the order of integrations in the two-dimensional integral, we 
get
\begin{equation}
\zeta(n)=\int\limits_0^1\frac{Li_{n-2}(x_2)}{x_2}\,dx_2\int\limits_{x_2}^1
\frac{dx_1}{x_1}=-\int\limits_0^1\frac{\ln{(x_2)}Li_{n-2}(x_2)}{x_2}\,dx_2.
\label{eq47}
\end{equation}
Now we can repeatedly apply the recurrence relation (\ref{eq24}), along 
with $Li_1(x)=-\ln{(1-x)}$ at the last step, and transform (\ref{eq47}) into
$$\zeta(n)=\int\limits_0^1\frac{\ln{x_1}}{x_1}\,dx_1\int\limits_0^{x_1}
\frac{dx_2}{x_2}\ldots\int\limits_0^{x_{n-4}}\frac{dx_{n-3}}{x_{n-3}}
\int\limits_0^{x_{n-3}}\frac{\ln{(1-x_{n-2})}}{x_{n-2}}\,dx_{n-2},$$
which after rescaling
$$x_2=x_1y_2,\;\;x_3=x_2y_3=x_1y_2y_3,\ldots,\;\; x_{n-2}=x_{n-3}y_{n-2}=
x_1y_2\cdots y_{n-2},$$
takes the form
\begin{equation}
\zeta(n)=\int\limits_0^1\frac{\ln{x_1}}{x_1}\,dx_1\int\limits_0^1\frac{dy_2}
{y_2}\ldots\int\limits_0^1\frac{dy_{n-3}}{y_{n-3}}\int\limits_0^1
\frac{\ln{(1-x_1y_2\cdots y_{n-2})}}{y_{n-2}}\,dy_{n-2}.
\label{eq48}
\end{equation}
Then the relation 
$$\int\limits_0^1\frac{dy_{n-1}}{1-x_1y_2\cdots y_{n-1}}=-
\frac{\ln{(1-x_1y_2\cdots y_{n-2})}}
{y_1y_2\cdots y_{n-2}}$$
shows that (\ref{eq48}) is equivalent to the $(n-1)$-dimensional integral
\begin{equation}
\zeta(n)=-\idotsint\limits_{\Box_{n-1}}\frac{\ln{x_1}}{1-x_1\cdots x_{n-1}}
\,dx_1\cdots dx_{n-1}.
\label{eq49}
\end{equation}
As in the previous case, (\ref{eq49}) can be further transformed into
$$\zeta(n)=-\frac{2^n}{2^n-1}\idotsint\limits_{\Box_{n-1}}\frac{\ln{x_1}}
{1-x^2_1\cdots x^2_{n-1}}\,dx_1\cdots dx_{n-1},$$
or, in the more symmetrical way,
\begin{equation}
\zeta(n)=-\frac{2^n}{2^n-1}\,\frac{1}{n-1}\,\idotsint\limits_{\Box_{n-1}}
\frac{\ln{(x_1\cdots x_{n-1})}}{1-x^2_1\cdots x^2_{n-1}}\,dx_1\cdots dx_{n-1}.
\label{eq50}
\end{equation}
Let us now assume that $n$ is odd and apply the Beukers-Kolk-Calabi change of 
variables (\ref{eq28}) to the integral (\ref{eq50}). We get
$$\zeta(2n+1)=-\frac{1}{2n}\,\frac{2^{2n+1}}{2^{2n+1}-1}\,
\idotsint\limits_{\Delta_{2n}}\ln{[\tan{(u_1)}\cdots \tan{(u_{2n})}]}\,
du_1\cdots du_{2n},$$
which is the same as
$$\zeta(2n+1)=-\frac{2^{2n+1}}{2^{2n+1}-1}
\idotsint\limits_{\Delta_{2n}}\ln{[\tan{(u_1)}]}\,du_1\cdots du_{2n}.$$
By rescaling variables, we can go from the polytope $\Delta_{2n}$ to the 
polytope $\delta_{2n}$ in this $2n$-dimensional integral and get
\begin{equation}
\zeta(2n+1)=-\frac{2^{2n+1}}{2^{2n+1}-1}\,\left(\frac{\pi}{2}\right)^{2n}
\idotsint\limits_{\delta_{2n}}\ln{\left[\tan{\left(u_1\frac{\pi}{2}\right)}
\right]}\,du_1\cdots du_{2n}.
\label{eq51}
\end{equation}
Using the kernel $K_{2n}(u,v)$, we can reduce the evaluation of (\ref{eq51})
to the evaluation of the following one-dimensional integral:
\begin{equation}
\zeta(2n+1)=-\frac{2\pi^{2n}}{2^{2n+1}-1}\int\limits_0^1
\ln{\left[\tan{\left(\frac{\pi}{2}u\right)}\right]}\,K_{2n}(u,u)\,du.
\label{eq52}
\end{equation}
But
$$\ln{\left[\tan{\left(\frac{\pi}{2}(1-u)\right)}\right]}=
\ln{\left[\cot{\left(\frac{\pi}{2}u\right)}\right]}=
-\ln{\left[\tan{\left(\frac{\pi}{2}u\right)}\right]},$$
which enables to rewrite (\ref{eq52}) as
\begin{equation}
\zeta(2n+1)=-\frac{\pi^{2n}}{2^{2n+1}-1}\int\limits_0^1
\ln{\left[\tan{\left(\frac{\pi}{2}u\right)}\right]}\,\left [ 
K_{2n}(u,u)- K_{2n}(1-u,1-u)\right ]\,du.
\label{eq53}
\end{equation}
However, from (\ref{eq42}) and (\ref{eq43}) we have (recall that 
$E_{2n-1}(1-u)=-E_{2n-1}(u)$)
$$K_{2n}(u,u)- K_{2n}(1-u,1-u)=(-1)^n\,\frac{2^{2n-1}}{(2n)!}\,
\frac{d}{du}E_{2n}(u),$$
and the straightforward integration by parts in (\ref{eq53}) yields finally 
the result
\begin{equation}
\zeta(2n+1)=\frac{(-1)^n\,\pi^{2n+1}}{4\,[1-2^{-(2n+1)}]\,(2n)!}\,
\int\limits_0^1\frac{E_{2n}(u)}{\sin{(\pi\,u)}}\,du.
\label{eq54}
\end{equation}
This is exactly the integral representation for $\zeta(2n+1)$ found in 
\cite{10}. Our earlier result (\ref{eq22}) for $\zeta(3)$ is just a 
special case of this more general formula.

\section{concluding remarks: $\zeta(2)$ and amoebas}
It remains to clarify the origin of the  highly non-trivial and miraculous
Beukers-Kolk-Calabi change of variables (\ref{eq28}). Maybe an interesting 
observation due to Passare \cite{7} that $\zeta(2)$ is related to the amoeba 
of the polynomial $1-z_1-z_2$ gives a clue.

Amoebas are fascinating objects in complex geometry \cite{18,19}. They are 
defined as follows \cite{20}. For a Laurent polynomial $P(z_1,\ldots,z_n)$,
let $Z_P$ denote the zero locus of $P(z_1,\ldots,z_n)$ in 
$(\mathbb{C}\backslash \{0\})^n$ defined by $P(z_1,\ldots,z_n)=0$. The amoeba 
${ A}(P)$ of the Laurent polynomial $P(z_1,\ldots,z_n)$ is the image of 
the complex hypersurface $Z_P$ under the map 
$$\mathrm{Log}:\; (\mathbb{C}\backslash \{0\})^n\to \mathbb{R}^n$$ 
defined through 
$$(z_1,\ldots,z_n) \to (\ln{|z_1|},\ldots,\ln{|z_n|}).$$

Let us find the amoeba of the following Laurent polynomial
\begin{equation}
P(z_1,z_2)=z_1-z_1^{-1}-i\left (z_2-z_2^{-1}\right ).
\label{eq55}
\end{equation}
Taking
$$z_1=e^u\,e^{i\phi_u},\;\;\; z_2=e^v\,e^{-i\phi_v},$$
we find that the zero locus of the polynomial (\ref{eq55}) is determined by 
conditions
$$\cos{\phi_u}\,\sinh{u}=\sin{\phi_v}\,\cosh{v},\;\;
\sin{\phi_u}\,\cosh{u}=\cos{\phi_v}\,\sinh{v}.$$
If we rewrite these conditions as follows
\begin{equation}
x=\frac{\sinh{v}}{\cosh{u}}=\frac{\sin{\phi_u}}{\cos{\phi_v}},\;\;\;
y=\frac{\sinh{u}}{\cosh{v}}=\frac{\sin{\phi_v}}{\cos{\phi_u}},
\label{eq56}
\end{equation}
we immediately recognize the Beukers-Kolk-Calabi substitution (\ref{eq8}) and
its hyperbolic version with the only difference that in (\ref{eq8}) we had 
$0\le x,y \le 1$. However, from (\ref{eq56}) we get
\begin{equation}
\cos^2{\phi_u}=\frac{1-x^2}{1-x^2y^2},\;\;\; 
\cos^2{\phi_v}=\frac{1-y^2}{1-x^2y^2},
\label{eq57}
\end{equation}
and
\begin{equation}
\cosh^2{u}=\frac{1+y^2}{1-x^2y^2},\;\;\; 
\cosh^2{v}=\frac{1+x^2}{1-x^2y^2}.
\label{eq58}
\end{equation}
It is clear from (\ref{eq57}) and (\ref{eq58}) that we must have
$$x^2\le 1,\;\;\; y^2\le 1.$$
Therefore, the amoeba ${ A}(P)$ is given by relations 
\begin{equation}
{ A}(P)=\left \{(u,v):\; -1\le\frac{\sinh{u}}{\cosh{v}}\le 1,\;\;
-1\le\frac{\sinh{v}}{\cosh{u}}\le 1\right \},
\label{eq59}
\end{equation}
and the  hyperbolic version of the Beukers-Kolk-Calabi change of variables
(\ref{eq8}) transforms the unit square $\Box$ into one-quarter of the 
amoeba (\ref{eq59}). Then the analog of (\ref{eq9}) indicates that $\zeta(2)$
equals one-third of the area of this amoeba.

As we see, the  hyperbolic version of the Beukers-Kolk-Calabi change of 
variables seems more fundamental and arises quite naturally in the context
of the amoeba (\ref{eq59}). Trigonometric version of it then is just an 
area-preserving transition from the ``radial'' coordinates $(u,v)$ to the 
``angular'' ones $(\phi_u,\phi_v)$.

Another amoeba related to $\zeta(2)$ was found in \cite{7}. Although the
corresponding amoeba ${ A}(1-z_1-z_2)$ looks different from the amoeba
(\ref{eq59}), they do have the same area. The trigonometric change of variables
used by Passare in \cite{7} is also different from (\ref{eq8}) but also leads
to simple calculation of the area of ${ A}(1-z_1-z_2)$ and hence 
$\zeta(2)$. Of course it will be very interesting to generalize this 
mysterious relations between $\zeta(n)$ and amoebas for $n>2$ and finally
disentangle the mystery. I'm afraid, however, that this game is already not 
for kids under fifteen.

\section*{Acknowledgments}
The author is grateful to Professor Noam D.~Elkies for his helpful comments
and suggestions. The work is supported by the Ministry of Education and 
Science of the Russian Federation and in part by Russian Federation President
Grant for the support of scientific schools NSh-2479.2014.2.

\end{document}